\documentclass[11pt]{paper}

\usepackage{amssymb,amsmath,enumerate,theorem,epsfig,rotating,graphics,changebar,eepic}
\usepackage{amsfonts}
\usepackage{a4,latexsym,parskip}
\usepackage{hyperref}
\hypersetup{pdfauthor=MS} \hypersetup{pdftitle=Fields}

\newcommand{\PP}{\mathbb{P}}

\newcommand{\C} [1][]{\mathbb{C}^{#1}}
\newcommand{\Q} [1] []{\mathbb{Q}_{#1}}
\newcommand{\N} [1][] {\mathbb{N}_{#1}}

\newcommand{\F}{\mathbb{F}}
\newcommand{\Z}{\mathbb{Z}}
\newcommand{\p}{\mathfrak{p}}
\newcommand{\m}{\mathfrak{m}}
\newcommand{\OO}{\mathcal{O}}

\newcommand{\qed}{\hfill \ensuremath{\Box}}
\newcommand{\NS}{\mathop{\rm NS}\nolimits}

\theoremstyle{break} \newtheorem{Theorem}{Theorem}
\newtheorem{Proposition}[Theorem]{Proposition}
\newtheorem{Lemma}[Theorem]{Lemma}

\newtheorem{Corollary}[Theorem]{Corollary}

\newtheorem{Remark}[Theorem]{Remark}

\newtheorem{Conjecture}[Theorem]{Conjecture}
\newtheorem{Observation}[Theorem]{Observation}

\newtheorem{Example}[Theorem]{Example}

\begin{document}
\setlength{\unitlength}{1cm}

\title{K3 surfaces with Picard rank 20} 

\author{Matthias Sch\"utt}


\date{\today}
\maketitle


\abstract{We determine all complex K3 surfaces with Picard rank 20 over $\Q$. Here the N\'eron-Severi group has rank 20 and is generated by divisors which are defined over $\Q$. Our proof uses modularity, the Artin-Tate conjecture and class group theory. With different techniques, the result has been established by Elkies to show that Mordell-Weil rank 18 over $\Q$ is impossible for an elliptic K3 surface. We then apply our methods to general singular K3 surfaces, i.e.~with N\'eron-Severi group of rank 20, but not necessarily generated by divisors over $\Q$.}




\section{Introduction}

Complex K3 surfaces of geometric Picard number 20 are called \emph{singular} since they involve no moduli. They share many properties with elliptic curves with complex multiplication (CM). For instance, they can always be defined over some number field. Moreover, over some finite extension of the number field, the $L$-series is given in terms of Hecke characters (cf.~Theorem~\ref{Thm:SI}). 

For singular K3 surfaces over $\Q$, Livn\'e proved the motivic modularity in \cite{L}. However, this definition does not require that the N\'eron-Severi group is generated by divisors which are defined over $\Q$. We refer to this particular property as ''Picard rank 20 over $\Q$".

The motivation to study such K3 surfaces was the following: In \cite{Shioda-20}, Shioda raised the question whether it was possible for an elliptic K3 surface to have Mordell-Weil rank 18 over $\Q$. One  way to disprove this would have been to show that in general, K3 surfaces with Picard rank 20 over $\Q$ do not exist.

However, it turned out that there are such examples (see Examples \ref{Ex:19}, \ref{Ex:7}). Recently Elkies determined all these surfaces in terms of their transcendental lattices:

\begin{Theorem}[Elkies {\cite{Elkies}}]\label{thm}
Let $X$ be a K3 surface with Picard rank 20 over $\Q$. Then the transcendental lattice $T(X)$ is primitive of class number one.
\end{Theorem}

Using sphere packings and gluing up to a Niemeier lattice, Elkies concluded that Mordell-Weil rank 18 over $\Q$ is impossible for an elliptic K3 surface.

Conversely, let $T(X)$ be primitive of class number one. Then the singular K3 surface $X$ with transcendental lattice $T(X)$ has a model with Picard rank 20 over $\Q$ (cf.~section~\ref{s:existence}).

In this paper, we present an alternative proof of  Theorem~\ref{thm} that we hope to be of independent interest. Our proof uses the following ingredients: modularity plus the classification of CM-forms in \cite{S-CM}; reduction and the Artin-Tate conjecture at split primes; class group theory.

We then generalise our techniques to all singular K3 surfaces. We deduce the following obstruction to the field of definition:

\begin{Theorem}\label{thm:2}
Let $L$ be a number field and $X$ a K3 surface of Picard rank 20 over $L$. Denote the discriminant of $X$ by $d<0$. Then $L(\sqrt{d})$ contains the ring class field $H(d)$.
\end{Theorem}

This result enables us to give a direct proof of \v Safarevi\v c' finiteness theorem for singular K3 surfaces (Theorem~\ref{Thm:oo}). It is the only known obstruction for the field of definition of a singular K3 surface other than the result on the genus of $T(X)$ in \cite{S-fields} (cf.~(\ref{eq:Cl}) and Lemma~\ref{Lem:2}). In a private correspondence, Elkies informed the author that his  proof for  Theorem~\ref{thm} also generalises to  Theorem~\ref{thm:2}.

The paper is organised as follows: The next two sections recall the relevant facts about singular K3 surfaces and modularity. In section \ref{s:Ex} we give two explicit examples of K3 surfaces of Picard rank 20 over $\Q$. Section \ref{s:AT} introduces the main techniques to be used, in particular the Artin-Tate conjecture. The proof of  Theorem~\ref{thm} is presented in the sections \ref{s:K} to \ref{s:prim}. The converse statement of  Theorem~\ref{thm} is covered in section \ref{s:existence}. We continue with the classification of K3 surfaces of Picard rank 20 over $\Q$ up to $\Q$-isomorphism. Section \ref{s:quad} generalises  Theorem~\ref{thm} to K3 surfaces with Picard rank 20 over a quadratic extension of $\Q$. The paper concludes with the proof of the general case of  Theorem~\ref{thm:2}.

\section{Singular K3 surfaces}
\label{s:sing}

The main invariant of a singular K3 surface $X$ is its \emph{transcendental lattice} $T(X)$. Here we consider the N\'eron-Severi group $\NS(X)$ of divisors up to algebraic equivalence as a lattice in $H^2(X, \Z)$ with cup-product. Then the transcendental lattice is the orthogonal complement
\[
T(X)=\NS(X)^\bot\subset H^2(X, \Z).
\]
The following classification was first stated by Pjatecki\u\i -\v Sapiro and  \v Safarevi\v c \cite{PSS}. The proof was completed by Shioda and Inose \cite{SI}:

\begin{Theorem}[Pjatecki\u\i -\v Sapiro - \v Safarevi\v c, Shioda - Inose]\label{Thm:Torelli}
The map $X\mapsto T(X)$ gives a bijection
\[
\{\text{Singular K3 surfaces}\}_{/\cong} \stackrel{1:1}{\longleftrightarrow} \{\text{positive-definite oriented even lattices of rank two}\}_{/\cong}.
\]
\end{Theorem}
The injectivity of this map follows from the Torelli theorem for singular K3 surfaces \cite{PSS}. For the surjectivity, Shioda-Inose exhibited an explicit construction involving isogenous CM-elliptic curves $E, E'$\cite{SI}. This is often referred to as Shioda-Inose structure:

$$
\begin{array}{ccccc}
E\times E' &&&& X\\
& \searrow && \swarrow &\\
&& \text{Km}(E\times E') &&
\end{array}
$$

Here both rational maps are 2:1, and $T(X)\cong T(E\times E')$. Shioda-Inose exhibited the rational map $X\dasharrow \text{Km}(E\times E')$ through base change of elliptic fibrations. Explicit equations were subsequently given by Inose in \cite{Inose}. In \cite{S-fields}, Inose's results were improved to derive a model over the ring class field $H(d)$ associated to the discriminant $d=\text{disc}(T(X))$ of the transcendental lattice (cf.~Lemma~\ref{Lem:def}). Over some extension, one can moreover determine the $\zeta$-function of $X$ (Theorem~\ref{Thm:SI}).

The set of singular K3 surfaces over $\Q$ (up to $\C$-isomorphism) is finite by a result of \v Safarevi\v c \cite{Shafa} (cf.~Theorem~\ref{Thm:oo}). However, there is only one effective obstruction known for a singular K3 surface $X$ to be defined over $\Q$: By \cite{S-fields}, the genus of $T(X)$ has to consist of a single class. (In \cite{Shimada}, Shimada proved this first for the case of fundamental discriminant $d$.) In other words, we require that its class group is only two-torsion:
\begin{eqnarray}\label{eq:Cl}
Cl(T(X))\cong (\Z/2)^g.
\end{eqnarray}
The general case will be treated in section~\ref{s:gen}. There we will also provide a formulation in terms of fields of definition (Lemma~\ref{Lem:2}).

The only drawback of relation (\ref{eq:Cl}) is that the class group $Cl(T(X))$ does not recognise whether $T(X)$ is primitive. We know 101 discriminants $d<0$ such that the class group $Cl(d)$ is only two-torsion. By a result of Weinberger \cite{Wb} there is at most one more such $d$, and in fact none under some condition on Siegel-Landau zeroes (which would follow from GRH). However, so far we lacked bounds for the degree of primitivity of $T(X)$. For Picard rank 20 over $\Q$, primitivity is part of  Theorem~\ref{thm}. For the general case, bounds for the degree of primitivity follow from  Theorem~\ref{thm:2} (cf.~section~\ref{s:gen}).

\section{Modularity of singular K3 surfaces over $\Q$}

We shall now see that condition (\ref{eq:Cl}) can also be understood in terms of modularity. Here the modular motive is the compatible system of Galois representations attached to the transcendental lattice. Over some extension, this motive is related to Hecke characters by  Theorem~\ref{Thm:SI}.

Throughout the paper, we fix the imaginary quadratic field $K=\Q(\sqrt{d})$ where $d<0$ is the discriminant of $X$. Write $d_K$ for the discriminant of $K$. Hence $d=N^2 d_K$.

\begin{Theorem}[Livn\'e {\cite{L}}]\label{Thm:Livne}
Every singular K3 surface $X$ over $\Q$ is modular. The $L$-series of the transcendental lattice $T(X)$ is the Mellin transform of a Hecke eigenform of weight 3 with CM by $K$.
\end{Theorem}

By a result of Ribet \cite{R}, CM-newforms are associated to Hecke characters. Essentially, a Hecke character $\psi$ of $K$ is given by its conductor $\m$, an ideal in the ring of integers $\OO_K$, and by its $\infty$-type $l$. Then $\psi$ satisfies
\[
\psi(\alpha\OO_K)=\alpha^l \;\;\;\forall\; \alpha\equiv 1\mod\m.
\]
Let $\mbox{N}_{K/\Q}$ denote the norm of $K/\Q$. The sum over all ideals $\mathfrak a$ of  $\OO_K$ that are relatively prime to $\m$ gives the $L$-function of $\psi$:
\[
L(\psi,s) = \sum_{\mathfrak a} \psi(\mathfrak a) \mbox{N}_{K/\Q}(\mathfrak a)^s.
\]
Through the inverse Mellin transform, $L(\psi,s)$ defines a newform of weight $k=l+1$ and level $|\mbox{N}_{K/\Q}(\m) d_K|$.
For the weight of the corresponding newform to be $3$, the Hecke character thus has to have $\infty$-type $2$. Moreover, we require the newform to have Fourier coefficients in $\Z$. This is possible if and only if the class group of $K$ consists only of two-torsion (cf.~Theorem~\ref{Thm:CM}). This condition is necessarily satisfied if (\ref{eq:Cl}) holds.

\begin{Example}\label{Ex:Hecke}
Let $K$ such that $Cl(K)\cong(\Z/2)^g$ with $d_K\neq -3, -4$. Then we can define a Hecke character $\psi$ of $K$ with trivial conductor and $\infty$-type $2$ by setting
\[
\psi(\alpha\OO_K)=\alpha^2
\]
for every principal ideal in $\OO_K$ and choosing suitable values for a set of generators of $Cl(K)$. Explicitly, let throughout this paper 
\[
D=\begin{cases}
-d_K, & \text{if } 4\nmid d_K,\\
-\frac{d_K}4, & \text{if } 4\mid d_K.
\end{cases}
\]
Assume that $p=\p\bar\p$ splits in $K$. Since $d_K\neq -3, -4$, we can write $p^2$ uniquely as
\[
p^2=x^2+D y^2,  \;\;\; x, y\in \frac 12 \N.
\]
(Here $x, y\in\N$ unless $D=-d_K$.) Then $\psi(\p)=\pm (x\pm\sqrt{-D} y)$. For the corresponding newform $f=\sum a_n q^n$, we obtain 
\[
a_p=\pm 2x.
\]
Once a normalisation is fixed, $f$ has level $|d_K|$ and Fourier coefficients in $\Z$.
\end{Example}

The newforms arising from different normalisations (i.e.~different sign choices) are quadratic twists of each other. In general, consider a (quadratic) Dirichlet character $\chi$ and a newform $f=\sum a_n q^n$. Then we obtain the twisted Hecke eigenform
\begin{eqnarray}\label{eq:f-twist}
f\otimes\chi = \sum_n a_n \chi(n) q^n.
\end{eqnarray}
The classification in \cite{S-CM} says that the construction of Example~\ref{Ex:Hecke} produces all Hecke characters resp.~Hecke eigenforms with Fourier coefficients in $\Z$ after twisting: 

\begin{Theorem}[Sch\"utt]\label{Thm:CM}
Let $K$ be an imaginary quadratic field. Then all Hecke characters of $K$ with fixed $\infty$-type $l$ such that the corresponding newform $f$ has coefficients in $\Z$, are identified under twisting. Moreover, there is such a Hecke character if and only if $Cl(K)\subseteq (\Z/l)^g$ for some $g\in\N$.
\end{Theorem}

\begin{Remark}
If $d_K\neq -3, -4$, then we only have to consider quadratic twists. If $\chi$ is a quadratic Dirichlet character, then we twist the Hecke character by $\chi\circ$N$_{\Q}^K$. In terms of the associated newform $f$, this corresponds to the quadratic twist in (\ref{eq:f-twist}). For $d_K=-3, -4$, we also have to take cubic resp.~biquadratic twisting into account. All these twists have geometric equivalents. For instance, any quadratic Dirichlet character can be identified with a Legendre symbol $\left(\frac\delta\cdot\right)$ for some squarefree $\delta\in\Z$. Then consider an elliptic curve (or a general equation of this type)
\begin{eqnarray}\label{eq:twist}
E:\;\;\;y^2=g(x)\;\;\;\text{ and twist } \;\;\; E_\delta:\;\;\; \delta y^2=g(x).
\end{eqnarray}
For geometric equivalents of cubic and biquadratic twists, the reader is referred to Remark~\ref{Rem:twists}.
\end{Remark}

\section{K3 surfaces of Picard rank 20 over $\Q$: Examples}
\label{s:Ex}

In this section, we recall two of the most elementary examples of K3 surfaces of Picard rank 20 over $\Q$. Both use elliptic fibrations with section. For further examples, the reader is referred to section \ref{s:existence}.

\begin{Example}\label{Ex:19}
There is a unique complex elliptic K3 surface $X$ with a fibre of type $I_{19}$. The fibration can be defined over $\Q$. This follows from work of Hall \cite{Hall} and was studied in detail by Shioda in \cite{ShCR}. A simple explicit Weierstrass equation is derived in \cite{SS2}:
\begin{eqnarray}\label{eq:19}
X:\;\;\;y^2 = x^3 + (t^4+t^3+3 t^2+1)\, x^2 + 2 (t^3 + t^2 + 2 t)\, x + t^2 + t + 1.
\end{eqnarray}

Let $U$ denote the hyperbolic plane 
generated by a general fibre and the zero-section. 
It is immediate that the N\'eron-Severi lattice of $X$ (over $\bar\Q$) can be written as
\[
\NS(X)=U \oplus A_{18}(-1) = \begin{pmatrix} 0 & 1\\1 & 0\end{pmatrix} \oplus A_{18}(-1)
\]
In particular, $X$ is a singular K3 surface. The Picard rank of $X$ over $\Q$ is 20 if and only if the all components of the special fibre are defined over $\Q$ (i.e.~the special fibre has  split multiplicative reduction). This can be achieved by an appropriate twist as in (\ref{eq:twist}) and was first exhibited in \cite{ST}. The model in (\ref{eq:19}) has the fibre of type $I_{19}$ at $t=\infty$. The fibre is split multiplicative, so the Picard rank of the surface over $\Q$ is already 20. The corresponding Hecke eigenform has level $19$ by \cite{ST} (see \cite[Table 1]{S-CM}).
\end{Example}

The next example goes back to Tate \cite{Tate-AEC}. It has been studied very concretely by Hulek and Verrill in \cite{HV}. 

\begin{Example}\label{Ex:7}
Let $X$ denote the universal elliptic curve for $\Gamma_1(7)$. Since this group has genus 0, the base curve is $\PP^1$. One the other hand, the space of cusp forms $S_3(\Gamma_1(7))$ is one-dimensional, so $X$ has geometric  genus $p_g(X)=1$. It follows that $X$ is a K3 surface. By general theory, the elliptic surface $X$ has a model over $\Q$ with a section $P$ of order 7 also defined over $\Q$. Such a model was first given by Tate in \cite{Tate-AEC}:
\[
X:\;\;\; y^2 + (1+t-t^2) xy + (t^2-t^3) y = x^3 + (t^2-t^3) x^2.
\]
Here $P=(0,0)$ is a point of order 7. In the following, we shall employ an abstract approach to show that $X$ has Picard rank 20 over $\Q$.

The quotient of $X$ by translation by $P$ gives rise to another elliptic K3 surface after resolving singularities. Hence the configuration of singular fibres can only be [1,1,1,7,7,7]. In particular, $X$ is a singular K3 surface. We claim that the above model has Picard rank 20 over $\Q$. Equivalently, each reducible fibre is completely defined over $\Q$. To prove this, we show that $P$ meets each $I_7$ fibre in a different non-trivial component. 

We employ Shioda's theory of Mordell-Weil lattices and the height pairing \cite{ShMW}. As a torsion section, $P$ has height 0. Since $P$ does not meet the 0-section, we can compute the height directly as
\[
h(P)=4-(\text{correction terms for reducible fibres}).
\]
Here the correction terms are $\frac{n(7-n)}7$ according to the component $\Theta_n$ which $P$ meets (cyclically numbered so that the zero-section meets $\Theta_0$). The only way to obtain $h(P)=0$ is
\[
0=h(P)=4-\frac 67-\frac {10}7-\frac{12}7.
\]
Since $P$ intersects each $I_7$ fibre at a non-trivial component, these special fibres are split multiplicative. Moreover, as the components differ for each $I_7$ fibre, their cusps cannot be conjugate. Hence all fibre components are defined over $\Q$, and the claim follows.
\end{Example}

\begin{Remark}\label{Rem:EMS}
The same argument applies to other modular elliptic K3 surfaces, but not to all of them. For instance, the universal elliptic curve for $\Gamma(4)$ is a Kummer surface. Hence it cannot have Picard rank 20 over $\Q$ by the next remark. This argument will also be used in the proof of the primitivity of the transcendental lattice (Lemma~\ref{Lem:prim}). Alternatively, we could also argue with the Weil pairing. Since the Weil pairing has image $\mu_4$, the fourth roots of unity, we deduce that $MW(X/\Q)\subset \Z/4 \times \Z/2$. Then we apply the inverse argument of Example~\ref{Ex:7} to a 4-torsion section which is not defined over $\Q$. This implies that there are singular fibres which are not completely defined over $\Q$.
\end{Remark}

\begin{Remark}[Singular abelian surfaces]\label{Rem:ab}
The situation for abelian surfaces is different: Let $A$ be a singular complex abelian surface, i.e.~$\rho(A)=4$. Then $A\cong E\times E'$ for isogenous CM-elliptic curves $E, E'$ by a result of Shioda-Mitani \cite{SM}. However, as Shioda noted in \cite{Sh-Error}, Picard rank 4 over $\Q$ is impossible. This is a consequence of the cohomology structure of abelian varieties and carries over to Kummer surfaces (cf.~Remark~\ref{Rem:EMS} and Lemma~\ref{Lem:prim}).
\end{Remark}

\section{The Artin-Tate conjecture}
\label{s:AT}

Let $X$ be a K3 surface of Picard rank 20 over $\Q$. In order to prove  Theorem~\ref{thm}, we will consider the reductions of $X$ at the good primes $p$ that split in $K$ and apply the Artin-Tate Conjecture. 

Let $p$ be a prime of good reduction of $X$. Then the reduction morphism induces embeddings
\begin{eqnarray}\label{eq:NS-prim}
\NS(X/\Q) \hookrightarrow \NS(X/\F_p)\;\;\;\text{ and } \;\;\; \NS(X/\bar\Q) \hookrightarrow \NS(X/\bar\F_p)
\end{eqnarray}
which are isometries onto the image. For almost all $p$, these embeddings are primitive. This follows from Shimada's argumentation in \cite[\S 2.2]{Shimada}, since the proof for the case of supersingular reduction can be generalised directly. For the remainder of the paper, we will only consider \emph{good primes} where the reduction is good and the embeddings in (\ref{eq:NS-prim}) are primitive.

On $X/\F_p$ we have the Frobenius endomorphism Frob$_p$ raising coordinates to their $p$-th powers. We want to consider the induced action on cohomology. For this, we fix a prime $\ell\neq p$ and work with \'etale $\ell$-adic cohomology of the base change $\bar X=X_{\bar\F_p}$ to an algebraic closure $\bar\F_p$ of $\F_p$. Then we consider the induced map Frob$_p^*$ on $H_{\acute{e}t}^2(\bar X, {\Q}_\ell)$ and its reciprocal characteristic polynomial
\[
P(X/\F_p, T)=\det(1-\text{Frob}_p^*\,T; H_{\acute{e}t}^2(\bar X, {\Q}_\ell)).
\]
Frob$_p^*$ acts through a  permutation on the algebraic cycles in $H_{\acute{e}t}^2(\bar X, {\Q}_\ell)$. More precisely, it operates as multiplication by $p$ on $\NS(X/\F_p)$ and in particular on the image of $\NS(X/\Q)$ under the primitive embedding (\ref{eq:NS-prim}). In the present case, $X$ has Picard rank 20 over $\Q$ and discriminant $d$. Let $f=\sum a_nq^n$ denote the associated newform by  Theorem~\ref{Thm:Livne}. Then
\begin{eqnarray}\label{eq:P_2}
P(X/\F_p, T)=(1-p\,T)^{20} \left(1- a_p\,T+\left(\frac dp\right)p^2\,T^2\right).
\end{eqnarray}
The Tate Conjecture \cite{Tate-C} relates the shape of the zeroes of $P(X/\F_p,T)$ to the Picard number: Conjecturally for any smooth projective surface $X$ over $\F_p$, it predicts
\begin{eqnarray*}
\rho(X/\F_p) & = & \# \{\text{zeroes } T=\frac 1p \text{ of } P(X/\F_p, T)\},\\
\rho(X/\bar\F_p) & = & \# \{\text{zeroes } T=\zeta\frac 1p \text{ of } P(X/\F_p, T) \text{ where $\zeta$ is a root of unity}\}.
\end{eqnarray*}
Here we count the zeroes with multiplicities. Since Frob$_p$ acts as multiplication by $p$ on $NS(X/\F_p)$, we always have $\leq$ in the above equations. For instance, the Tate conjecture is known for elliptic K3 surfaces \cite{ASD}. By \cite{Milne} (cf.~\cite[p.~25]{Milne-A} for characteristic two), it is equivalent to the Artin-Tate Conjecture:

\begin{Conjecture}[Artin-Tate {\cite{Artin-Tate}}]
Let $X/\F_p$ be a smooth projective surface. Let $\alpha(X)=\chi(X) -1\, +$ dim Pic Var$(X)$. Then
\begin{eqnarray}\label{eq:ATC}
\dfrac{P(X/\F_p, T)}{(1-p\,T)^{\rho(X/\F_p)}}\Big |_{T=\frac 1p} =  \dfrac{|Br(X/\F_p)|\;|\text{discr}(\NS(X/\F_p))|}{p^{\alpha(X)}\; |\NS(X/\F_p)_\text{tor}|^2}
\end{eqnarray}
\end{Conjecture}

\begin{Remark}\label{Rem:Brauer}
By \cite{Brauer}, $|Br(X/\F_p)|$ is always a square. For K3 surfaces, $\alpha(X)=1$ and the N\'eron-Severi group is torsion-free, since numerical and algebraic equivalence coincide. Hence (\ref{eq:ATC}) simplifies to 
\begin{eqnarray}\label{eq:AT}
\dfrac{P(X/\F_p, T)}{(1-p\,T)^{\rho(X/\F_p)}}\Big |_{T=\frac 1p} =  \frac 1p\; |Br(X/\F_p)|\;|\text{discr}(\NS(X/\F_p))|.\end{eqnarray}
\end{Remark}

We shall now specialise to the situation where $X$ is a K3 surface with Picard rank 20 over $\Q$ and $p$ is a good split prime. The Fourier coefficient $a_p$ can be computed in terms of Example~\ref{Ex:Hecke}. In particular, it is never a multiple of $p$. Hence the zero $T=\frac 1p$ of $P(X/\F_p, T)$ has multiplicity exactly 20, and there is no further zero $T=\zeta\frac 1p$. It follows that $\rho(X/\F_p)=\rho(X/\bar\F_p)=20$. In particular, the Tate conjecture  holds for $X$ over $\F_p$. From (\ref{eq:NS-prim}) we deduce
\begin{eqnarray*}
\NS(X/\bar\Q) = \NS(X/\Q) = \NS(X/\F_p) = \NS(X/\bar\F_p)
\end{eqnarray*}
and thus
\begin{eqnarray*}
\text{discr}(\NS(X/\Q))=\text{discr}(\NS(X/\F_p))=d=N^2\,d_K.
\end{eqnarray*}
Hence the Artin-Tate Conjecture for $X/\F_p$ (\ref{eq:AT}) gives with $M^2=|Br(X/\F_p)|$ 
\begin{eqnarray}\label{eq:d}
2\,p-a_p = M^2 |d| = (MN)^2 |d_K|.
\end{eqnarray}

The proof of  Theorem~\ref{thm} now proceeds in three steps:
\begin{enumerate}[1.]
\item The imaginary quadratic field $K$ has class number one (Corollary~\ref{Cor:K}).
\item The discriminant  $d$ has class number one (Corollary~\ref{Cor:d}).
\item The transcendental lattice $T(X)$ is primitive (Lemma~\ref{Lem:prim}).
\end{enumerate}
As a by-product, we will also determine the possible shapes of the associated newform $f$ (Lemma~\ref{Lem:f}).

\section{Class number of $K$}
\label{s:K}

In this section, we will prove that $K$ has class number one. We achieve this through the following Proposition:

\begin{Proposition}\label{Prop:K}
Let $p$ split in $K$ and $a_p\in\Z$ the coefficient of a newform of weight $3$ with CM by $K$. Then (\ref{eq:d}) implies that $p$ splits into principal ideals in $K$.
\end{Proposition}

\emph{Proof:} By Example~\ref{Ex:Hecke}, we can write $a_p=2z$ with $z=\pm x\in\frac 12\Z$. By (\ref{eq:d}), we have
\begin{eqnarray}\label{eq:z}
p-z & = & \frac{m^2 D}2
\end{eqnarray}
for some $m\in\N$. On the other hand, $p^2=z^2+D y^2$ for some $y\in\frac 12\N$ by assumption, i.e.
\begin{eqnarray}\label{eq:y}
p^2 - z^2 &  = & D y^2.
\end{eqnarray}
Dividing (\ref{eq:y}) by (\ref{eq:z}), we obtain
\begin{eqnarray}\label{eq:p}
p+z & = & 2\left(\frac ym\right)^2.
\end{eqnarray}
Now we add (\ref{eq:z}) and (\ref{eq:p}) and  divide by two  to derive
\begin{eqnarray}\label{eq:pp}
p  & = & \left(\frac ym\right)^2 + D \left(\frac m2\right)^2.
\end{eqnarray}
Since $\frac m2\in\frac 12\N$, the same holds for $\frac ym$. We deduce that $p$ splits into principal ideals in $K$. \qed

\begin{Corollary}\label{Cor:K}
Let $X$ be a K3 surface of Picard rank 20 over $\Q$. Then its CM-field $K$ has class number one.
\end{Corollary}

\emph{Proof:} By the Artin-Tate conjecture, eq.~(\ref{eq:d}) holds at all but finitely many $p$ that split in $K$. By Proposition~\ref{Prop:K}, each of these $p$ splits into principal ideals in $K$. Hence $K$ has class number one. \qed

\section{Shape of $f$}

If $K$ has class number one, we can describe the CM-newforms of $K$ even more explicitly in terms of Example~\ref{Ex:Hecke}. Here we only have to take extra care of the special cases $d_K=-3, -4$ where $\OO_K\neq\{\pm 1\}$. For this purpose, let 
\[
D'=\begin{cases}
27, & \text{if $d_K=-3$,}\\
4, & \text{if $d_K=-4$,}\\
D, & \text{if $d_K\neq -3,-4$}.
\end{cases}
\]

\begin{Example}[Class number one]\label{Ex:h=1}
Let $K$ have class number one. Let $D'$ as above. If $p$ splits in $K$, then we rewrite (\ref{eq:pp}) uniquely as
\[
p=x^2+D'y^2\;\;\; x,y\in\frac 12\N.
\]
The corresponding Hecke character $\psi$ of $\infty$-type $2$ sends the prime ideal $(x+\sqrt{-D'}\,y)$ to its square. We obtain the newform $f_K$ of weight $3$ and level $D'$ from \cite[Tab.~1]{S-CM} with coefficients
\begin{eqnarray}\label{eq:f}
a_p=2(x^2-D'y^2).
\end{eqnarray}
\end{Example}

\begin{Lemma}\label{Lem:f}
Let $X$ be a K3 surface of Picard rank 20 over $\Q$. Let $f$ denote the associated newform.
\begin{enumerate}[(i)]
\item If $d_K\neq -3, -4$, then $f=f_K$. 
\item If $d_K=-4$, then $f$ is a quadratic twist of $f_K$. 
\item If $d_K=-3$, then $f$ is a cubic twist of $f_K$.  
\end{enumerate}
\end{Lemma}

\emph{Proof:} Assume that $d_K\neq -3, -4$. Let $p$ a split prime as in Example~\ref{Ex:h=1}. By  Theorem~\ref{Thm:CM}, $f$ has the coefficient
\begin{eqnarray}\label{eq:g}
a_p=\pm 2 (x^2-Dy^2).
\end{eqnarray}
Inserting into (\ref{eq:d}) gives
\begin{eqnarray}\label{eq:gg}
2 (x^2+Dy^2 \mp (x^2-Dy^2)) = m^2 D.
\end{eqnarray}
Since $d_K$ is not a square and neither is $D$, it follows that only the minus sign in (\ref{eq:gg}) is possible. I.e.~in (\ref{eq:g}), only the plus sign occurs. By definition $f=f_K$.

If $d_K=-4$ and $p=x^2+4y^2$, then
\[
a_p=\begin{cases}
\pm 2(x^2-4y^2),\\
\pm 8xy.
\end{cases}
\]
The second case occurs (at some split $p$) if and only if $f$ is a biquadratic twist of $f_K$. Only the first case is compatible with (\ref{eq:d}), since in the second case
\[
2p-a_p= 2(x^2+4y^2\mp 4xy)=2(x\mp 2y)^2\neq 4n^2.
\]
Hence $f$ is a quadratic twist of $f_K$.

A similar argument rules out quadratic and sextic twists of $f_K$ for $d_K=-3$: Here we can always write the coefficients of $f$ as
\[
a_p=\pm 2(x^2-3y^2)\;\;\;\text{ where non-uniquely}\;\;\; p=x^2+3y^2, \;\; x,y\in\frac 12 \N.
\]
By the argument of case (i), only the plus sign occurs. This implies that $f$ is a cubic twist of $f_K$.
 \qed

\section{Class number of $d$}
\label{s:d}

Let $X$ be a K3 surface of Picard rank 20 over $\Q$. Denote the associated newform by $f=\sum a_n q^n$. We can rephrase Lemma~\ref{Lem:f} and its proof as follows: At every good split prime $p$, we can write (non-uniquely if $D\neq D'$)
\begin{eqnarray}\label{eq:dd}
p=x_p^2+Dy_p^2\;\;\;\text{such that}\;\;\; a_p=2(x_p^2-Dy_p^2) \;\;\;\text{and}\;\;\; 4 D y_p^2 = M_p^2\, |d|.
\end{eqnarray}
By construction, we have either $d_K=-4D$ and $y\in\N$ or $d_K=-D$ and $y\in\frac 12 \N$. Recall that $d=N^2 d_K$ and $d_K$ has class number one by Corollary~\ref{Cor:K}. We want to find all $d$ which are compatible which Picard rank 20 over $\Q$. In other words, we search for all $N|M_p$ which are simultaneously possible in (\ref{eq:dd}) at all good split $p$.

\begin{Observation}\label{Obs}
Let $gcd$ denote the greatest common divisor in $\N$ if $d_K=-4D$, resp.~in $\frac 12 \N$ if $d_K=-D$. Let $y_p$ be given by (\ref{eq:dd}) at a good split prime $p$. Then
\[
N| \; gcd(y_p;  \,p \text{ good split prime for } X).
\]
\end{Observation}

Hence, if for instance there was one $y_p=1$ resp.~$y_p=\frac 12$ occurring, then $d=d_K$ (and $N=M_p=1$) would follow. However, this need not be the case in general. To see this, let the associated newform $f$ have level $27$. Then by construction $3|y_p$ for all split $p$. Hence at least $d=-3$ and $d=-27$ would be possible a priori. 

To bound $d$ (or $N$) in general, we need information on the greatest common divisor of the $y_p$. This divisibility problem translates into class group theory through representations of primes by quadratic forms:

\begin{Lemma}\label{Lem:r}
Let $d<0$ and $Q=\begin{pmatrix} 2 & b\\b & 2c\end{pmatrix}$ a quadratic form of discriminant $d$. Consider the primes $p$ represented by $Q$:
\begin{eqnarray}\label{eq:u,v}
p=u_p^2+b u_p v_p+c v_p^2\;\;\;\; u_p, v_p\in\N.
\end{eqnarray}
Then almost all $v_p$ are divisible by $r\in\N$ if and only if $h(d)=h(dr^2)$.
\end{Lemma}

\emph{Proof:} Note that $Q$ always represents the principal class in $Cl(d)$. Hence, if $h(d)=h(dr^2)$, then the quadratic form $Q_r=\begin{pmatrix} 2 & br\\br & 2cr^2\end{pmatrix}$ in $Cl(dr^2)$ represents the same primes as $Q$ (the principal ones). Thus $r|v_p$ for all these $p$.

Conversely, assume that $r|v_p$ for almost all $p$ represented by $Q$. Thus all these $p$ are represented by $Q_r$ as well. Since the split primes are equally distributed on the classes that represent them, we obtain $h(d)\geq h(dr^2)$. On the other hand, $h(d)\leq h(r^2 d)$ holds trivially. Hence the class numbers $h(d)$ and $h(dr^2)$ have to coincide. \qed

\begin{Corollary}\label{Cor:d}
Let $X$ be a K3 surface of Picard rank 20 over $\Q$. Then the transcendental lattice has discriminant $d$ of class number one.
\end{Corollary}

\emph{Proof:} By Corollary~\ref{Cor:K}, $d_K$ has class number one. Assume that $d\neq d_K$, i.e.~there is some $r$ dividing all $y_p$ in (\ref{eq:dd}). To apply Lemma~\ref{Lem:r}, we have to relate divisibility of $y_p$ and $v_p$. We consider the following quadratic forms:
\[
Q=\begin{pmatrix} 2 & 0\\0 & 2D\end{pmatrix}, \;\;\text{if $d_K$ is even,}\;\;\;\; Q=\begin{pmatrix} 2 & 1\\1 & \frac{D+1}2\end{pmatrix}, \;\;\text{if $d_K$ is odd.}
\]
If $d_K$ is even, then $d_K=-4D$ and $v_p=y_p\in\N$. Hence $h(d)=h(d_K)=1$ follows from Lemma~\ref{Lem:r}. If $d_K$ is odd, then we can rewrite (\ref{eq:u,v}) in half-integers:
\[
p=u_p^2+u_p v_p+\frac{D+1}4 v_p^2= \left(u_p+\frac {v_p}2\right)^2 + D \left(\frac {v_p}2\right)^2.
\]
Hence divisibility of $y_p$ in $\frac 12\N$ translates into divisibility of $v_p\in\N$ and vice versa. Again we deduce $h(d)=h(d_K)=1$ by Lemma~\ref{Lem:r}. \qed


\begin{Remark}
Let $X$ be a K3 surface of Picard rank 20 over $\Q$. If $d\neq d_K$, it is immediate from the above argument that the associated newform $f$ has a particular shape. For $d=-28$, this newform is uniquely determined with level $7$ by Lemma~\ref{Lem:f}. In the other three cases ($d=-12,-16,-27$), it is easily checked that the condition $r|y_p$ fixes a unique Hecke character. We find that $f$ is the unique newform of weight 3 and  level $|d|$.
\end{Remark}

\section{Primitivity of $T(X)$}
\label{s:prim}

We have seen that a K3 surface with Picard rank 20 over $\Q$ has discriminant of class number one. Hence there are a priori 17 possibilities for the transcendental lattice: 

\begin{enumerate}[(1)]
\item the 13 primitive lattices of class number one, corresponding to isomorphism classes of CM-elliptic curves over $\Q$ through the Shioda-Inose structure,
\item the four imprimitive lattices of discriminant $d=-12, -16, -27, -28$.
\end{enumerate}

In this section, we will rule out the second case. 

\begin{Lemma}\label{Lem:prim}
Let $X$ be a K3 surface of Picard rank 20 over $\Q$. Then $T(X)$ is primitive.
\end{Lemma}

\emph{Proof:} Assume that $T(X)$ is not primitive. By Corollary~\ref{Cor:d}, we are in case (2) above. We shall treat even and odd discriminants separately.

If $d$ is even in case (2) above, then the transcendental lattice $T(X)$ has intersection form $2 Q$ for $Q\in Cl(d')$ where $4d'=d$. It follows from \cite{SI} that $X$ is the Kummer surface of an abelian surface $A$ such that the transcendental lattice $T(A)$ has intersection form $Q$. By Remark~\ref{Rem:ab}, $\rho(A/\Q)<4$ and $\rho(X/\Q)\leq\rho(A/\Q)+16<20$.

If $d$ is odd, i.e.~$d=-27$, then we consider Inose's fibration on $X$ (cf.~\cite{Inose}, \cite{Sandwich}). In the present case, $K=\Q(\sqrt{-3})$, and $X$ arises from the Shioda-Inose construction (cf.~section~\ref{s:sing}) for the following elliptic curves:
\[
E \text{ with CM by } \OO_K\;\;\; \text{ and } \;\;\; E' \text{ with CM by } \Z+3\OO_K.
\] 
In particular, $j(E)=0$. It follows from \cite{Inose} that $X$ admits the isotrivial elliptic fibration
\[
X:\;\;\; y^2 = x^3 + t^5 (3\,t^2-2\cdot 11\cdot 23\, t+3).
\]
Here the singular fibres have type $II^*, II^*, II, II,$ and the Mordell-Weil group over $\bar\Q$ has rank two. The generic fibre has CM by $\OO_K$. Let $\omega$ denote a primitive third root of unity acting on $X$ via $x\mapsto \omega \,x$. If $P$ is a section of the elliptic surface, then so is $\omega^*\,P$. Since the singular fibres admit no non-trivial torsion sections, these sections are independent. Since this argumentation applies to any twist $Y$ of $X$, Gal$(\Q(\sqrt{-3})/\Q)$ always acts non-trivially on $MW(Y)$. Hence rk $MW(Y/\Q)<2$ and in particular $\rho(Y/\Q)<20$. \qed

\section{Existence of K3 surfaces of Picard rank 20 over $\Q$}
\label{s:existence}

There are 13 primitive lattices $T$ of class number one appearing in  Theorem~\ref{thm}. For each of them one can ask whether there is a K3 surface with Picard rank 20 over $\Q$ and this transcendental lattice. Elkies announced in \cite{Elkies} that this holds true for each $T$. It follows that for each of these surfaces, one such model is given by Inose's fibration for the CM-elliptic curve corresponding to $T$, as exhibited over $\Q$ in \cite{S-fields}. 

However, for Inose's fibration, the non-trivial sections are often not immediate. In the cases at hand, there are two fibres of type $II^*$ plus an additional reducible fibre of type $I_2$. Hence the Mordell-Weil rank is one. Elkies recently computed the Mordell-Weil generator of height $\frac{|d|}2$ explicitly for all these fibrations \cite{E-web}.

For the reader's convenience we include a list of different models of these K3 surfaces where the Picard rank 20 over $\Q$ becomes evident. These models are given in terms of elliptic fibrations with configuration of singular fibres and the abstract structure of the Mordell-Weil group. We also include a reference, but naturally the given models are far from unique. Other models may be found in \cite{Noam}, \cite{S-Rocky}, \cite{TY} for instance. Explanations follow the table.

$$
\begin{array}{|c|c|c|c|}
\hline
d & \text{configuration} & MW  & \text{reference}\\
\hline\hline
-3 & [1^3,3,12^*] & \Z/4  & \text{Lem}.~\ref{Lem:-3}\\
\hline
-4 & [0^*,III^*, III^*] & \Z/2 & \text{Lem.}~\ref{Lem:-4}\\ 
\hline
-7 & [1^3,7^3] & \Z/7  & \text{Ex.}~\ref{Ex:7}\\
\hline
-8 & [1,4,III^*,II^*] & \{0\}  & \cite[\S 7]{S-Rocky} \\
\hline
-11 & [1^3,11,II^*] & \{0\}  & \cite[(III.2)]{S-Diss}\\
\hline
-12 & [2,3,III^*, II^*] & \{0\}  & \cite[\S 7]{S-Rocky}\\
\hline
-16 & [2,8,1^*,1^*] & \Z/4  & \cite[\S 7]{S-Rocky}\\
\hline
-19 & [1^5,19] & \{0\}  & \text{Ex.}~\ref{Ex:19}\\
\hline
-27 & [1^4,2,9^2] & \Z + \Z/3 & \text{Ex.}~\ref{Ex:d=-27}\\
\hline
-28 & [1^6,6,12] & \Z^2  & \cite[\S 5]{Noam}\\
\hline
-43 & [1^6,6,12] & \Z^2  & \cite[\S 5]{Noam}\\
\hline
-67 & [1^3,4,7,II^*] & \Z  & \cite[\S 4]{Noam}\\
\hline
-163 & [1^6,6,12] & \Z^2  & \cite[\S 5]{Noam}\\
\hline
\end{array}
$$

For $\boldsymbol{d=-8, -12}$, it was shown in \cite[\S 7]{S-Rocky} that the named fibrations are defined over $
\Q$. To obtain Picard rank 20 over $\Q$, it suffices to apply a quadratic twist as in Example~\ref{Ex:19} such that the fibre of type $I_4$ resp.~$I_3$ becomes split-multiplicative.

For $\boldsymbol{d=-11}$, the following Weierstrass form was derived in \cite{S-Diss}:
\begin{eqnarray*}
y^2 = x^3+t^2(t^2+3t+1)\,x^2+t^4(2t+4)\,x+t^5(t+1).
\end{eqnarray*}
This fibration has a $II^*$ fibre at 0 and a split-multiplicative fibre of type $I_{11}$ at $\infty$. 

For $\boldsymbol{d=-16}$, we realise the surface as a  quadratic base change of the extremal rational elliptic surface with  configuration $[1,4,1^*]$. It has a rational 4-torsion section $P$ which meets the singular fibres $I_4$ at a near and $I_1^*$ at a far component (cf.~\cite{MP1}). This implies that all fibre components are defined over $\Q$. The same argumentation applies to the base changed surface. Here we choose the base change in such a way that the $I_1^*$ fibres sit above rational cusps.

\begin{Example}[Discriminant {$\boldsymbol{d=-27}$}]
\label{Ex:d=-27}
For this discriminant, we searched the one-dimensional family of elliptic K3 surfaces with the given configuration $[1^4, 2, 9^2]$ and a 3-torsion section for an appropriate specialisation. Using techniques from \cite{ES}, we found
\[
X:\;\;\; {y}^{2}+ 3 \left(2{t}^{2}+1 \right) xy+ \left(1-t^2 \right)^{3}y = {x}^{3}.
\]
This elliptic surface has 3-torsion sections with zero $x$-coordinate and an independent section $P$ over $\Q$ with $x$-coordinate $x(P)=(t-1)^3$ and height $h(P)=3/2$. The $I_9$ fibres are located at $t=\pm1$ and split-multiplicative. Hence $X$ has Picard rank 20 over $\Q$. Using the height pairing \cite{ShMW}, one can show that neither $P$ nor its translates by the torsion sections are 3-divisible. Hence $X$ has discriminant
\[
d=-h(P)\,\dfrac{\text{disc}(A_1)\,\text{disc}(A_8)^2}{|MW(X)|^2} = -27.
\]
\end{Example}

\section{Classification up to $\Q$-isomorphism}

So far, we have only considered K3 surfaces up to isomorphism over $\C$. Then a singular K3 surface $X$ is identified by its transcendental lattice $T(X)$ (Theorem~\ref{Thm:Torelli}). In this section, we shall answer the question which $\Q$-isomorphism classes have Picard rank 20 over $\Q$. This is closely related to the precise shape of the corresponding Hecke eigenform (cf.~Lemma \ref{Lem:f}). 

We will work with Inose's elliptic fibration. In this context, one should always have quadratic twisting as in (\ref{eq:twist}) in mind. This operation twists the modular forms. Notably it also twists sections and affects singular fibres of types $IV, IV^*, I_m^*, I_n (n>2)$. Our first result concerns the case $d\neq -3, -4$:

\begin{Proposition}\label{Prop:Q}
Let $0>d\neq -3, -4$ of class number one. Up to $\Q$-isomorphism, there is a unique K3 surface $X$ of discriminant $d$ and Picard rank 20 over $\Q$.
\end{Proposition}

\emph{Proof:} The existence was shown in the previous section. We work with Inose's fibration with reducible singular fibres $I_2, II^*, II^*$ and a section $P$ of height $h(P)=\frac{|d|}2$. Over $\C$, such a fibration is unique (cf.~\cite{Sandwich}). Over $\Q$, this only leaves quadratic twists (for $d\neq -3, -4$). But then the condition that the section $P$ is defined over $\Q$ distinguishes the unique twist with Picard rank 20 over $\Q$. \qed

\begin{Lemma}\label{Lem:-4}
Let $d=-4$.  Consider the extremal elliptic K3 surface
\begin{eqnarray}\label{eq:d=-4}
X:\;\;\;y^2 = x^3 - t^3(t-1)^2 x
\end{eqnarray}
with singular fibres $III^*$ at $0$ and $\infty$ and $I_0^*$ at $1$ and two-torsion section $(0,0)$. Then any K3 surface with discriminant $d$ and Picard rank 20 over $\Q$ is $\Q$-isomorphic to a quadratic twist of $X$.
\end{Lemma}

\emph{Proof:} The configuration determines a unique elliptic fibration over $\C$. Over $\Q$, we distinguish biquadratic twists
\[
X_\delta:\;\;\;y^2 = x^3 - \delta\, t^3(t-1)^2 x\;\;\;\;(\delta\in\Q^*).
\]
All fibre components are defined over $\Q$ with the possible exception of the simple components of  the $I_0^*$ fibre which do not meet the zero section. These components are endowed with the Galois action of the extension $\Q(x^3 - \delta\,x)/\Q$. Hence all components are defined over $\Q$ if and only if $\delta$ is a square in $\Q^*$. This corresponds to the quadratic twist of (\ref{eq:d=-4}) by $\sqrt{\delta}$ as in (\ref{eq:twist}). \qed

\begin{Lemma}\label{Lem:-3}
Let $d=-3$.  Consider Inose's fibration
\begin{eqnarray*}
X:\;\;\;y^2 = x^3 - t^5(t+1)^2
\end{eqnarray*}
with singular fibres $II^*$ at $0$ and $\infty$ and $IV$ at $-1$. Then any K3 surface with discriminant $d$ and Picard rank 20 over $\Q$ is $\Q$-isomorphic to a cubic twist of $X$.
\end{Lemma}

Different elliptic fibrations on this surface have been studied in \cite{S-K3}. We omit the proof which is analogous to the previous one.

\begin{Remark}\label{Rem:twists}
If $d=-3$ or $-4$, then there are infinitely many possible associated newforms by Lemma~\ref{Lem:f}. By the previous two lemmata, each of these twists (quadratic resp.~cubic) is associated to a unique K3 surface of Picard rank 20 over $\Q$. 
\end{Remark}


\section{K3 surfaces with Picard rank 20 over a quadratic extension}
\label{s:quad}

In the next section, we will apply our methods to fields of definition of general singular K3 surfaces and their N\'eron-Severi lattices. To give a flavor of the ideas involved, we first give a full treatment of K3 surfaces with Picard rank 20 over a quadratic extension of $\Q$. Throughout, we employ the same techniques and notation as above.

\begin{Proposition}\label{Prop:2}
Let $L$ be a quadratic extension of $\Q$ and $X$ be a K3 surface with Picard rank 20 over $L$. As before, let $T(X)$ denote the transcendental lattice, $d$ its discriminant and $K=\Q(\sqrt{d})$. Then:
\begin{enumerate}[(i)]
\item If $L=K$, then $d$ has class number one.
\item If $L\neq K$, then $d$ has class number one or two. In the latter case, the compositum $LK$ agrees with the ring class field $H(d)$.
\end{enumerate}
\end{Proposition}

\emph{Proof:} We shall consider all those primes $p$ that split in both $K$ and $L$. Let $\p|p$ in $L$. Then $\F_\p=\F_p$ and again $\rho(X/\F_p)=20$. As before we will apply the Artin-Tate Conjecture to the reduction of $X$ at $\p$. For this, we need the coefficient $a_\p$ of the characteristic polynomial of Frob$_\p$ as in (\ref{eq:P_2}). 
Even if $X$ is not defined over $\Q$, there still is a modularity result over some extension:

\begin{Theorem}[Shioda-Inose {\cite[Theorem~6]{SI}}]
\label{Thm:SI}
Upon increasing the base field, the $\zeta$-function of a singular K3 surface $X$ splits into one-dimensional factors. Then the $L$-function of the transcendental lattice factors as
\[
L(T(X),s) = L(\psi^2, s)\, L(\bar\psi^2, s)
\]
where $\psi$ is the Hecke character associated to an elliptic curve with CM in $K$. Here one can choose the elliptic curve $E$ identified with the transcendental lattice $T(S)$ under the map
\[
\begin{pmatrix}
2a & b\\ b & 2c
\end{pmatrix} \mapsto \tau=\dfrac{-b+\sqrt{b^2-4ac}}{2a} \mapsto E=\C/(\Z+\tau\Z).
\]
\end{Theorem}

Thanks to this result, we are able to derive the relevant properties of $a_\p$ to apply our previous techniques.
We will need the following lemma:

\begin{Lemma}\label{Lem:shape}
In the above notation, $a_\p\in K$. Moreover $a_\p$ takes the shape of  Example~\ref{Ex:Hecke} and $p$ splits into primes of order two in $Cl(K)$:
\[
p^2 = \alpha_\p \cdot \bar\alpha_\p = x^2 + D y^2,\;\;\;\; a_\p = \pm 2 x.
\]
\end{Lemma}

\begin{Remark}
\label{Rem:new}
By inspection,
Lemma \ref{Lem:shape} does not require $\NS(X)$ to be rational over $\F_\p$, but only $p$
to split in $K$ and $\p\mid p$.
In our special case where $\NS(X)$ is fully defined over $\F_\p$, 
the proof of Proposition \ref{Prop:2} will ultimately show that  $p$ splits into principal primes  in $Cl(d)$.
\end{Remark}

\emph{Proof of Lemma \ref{Lem:shape}:} From the Weil conjectures we know that 
\begin{eqnarray}\label{eq:alpha_p}
a_\p=\alpha_\p+\bar\alpha_\p\;\;\;\text{with}\;\; |\alpha_\p|=p.
\end{eqnarray}
Here $\alpha_\p, \bar\alpha_\p$ are algebraic integers, complex conjugate in an imaginary quadratic extension of $\Q$ since $a_\p\in\Z$. We have to show that this quadratic field is $K$. In the present situation we know the $\zeta$-function of $X$ over some extension of $L$ by Theorem \ref{Thm:SI}.
As a result of increasing the ground field, the eigenvalues $\psi(\mathfrak{P})^2, \overline{\psi(\mathfrak{P})}^2$ of Frobenius at a prime $\mathfrak{P}$ above $\p$ agree with some power of $\alpha_\p, \bar\alpha_\p$. Since $\psi(\mathfrak{P})\in K\setminus\Q$ and $\alpha_\p$ is quadratic over $\Q$, this implies that $\alpha_\p\in K$. It follows that $a_\p$ has exactly the same shape as $a_p$ in Example~\ref{Ex:Hecke}. In fact, we deduce from (\ref{eq:alpha_p}) that
\[
p^2 = \alpha_\p \cdot \bar\alpha_\p = x^2 + D y^2,\;\;\;\text{ where } a_\p=\alpha_\p+\bar\alpha_\p=\pm 2x.
\]
This is to say that the prime factors of $p$ in $K$ become principal upon squaring. \qed

Thanks to Lemma \ref{Lem:shape} we can continue exactly along the lines of the previous sections to complete the proof of Proposition~\ref{Prop:2}. We distinguish two cases:

If $L=K$, then at every good split prime $\p$ in $K$, we have $\rho(X/\F_\p)=20$. Hence the arguments from the previous sections carry over except for Lemma~\ref{Lem:prim}. I.e., $d$ has class number one, but imprimitive $T(X)$ occurs.

If $L\neq K$, then Proposition~\ref{Prop:K} tells us that all the primes that split in both $K$ and $L$ are principal. Hence $K$ has class number one or two. By the argumentation of section \ref{s:d}, all these $p$ are principal in $Cl(d)$ as well (as mentioned in~Remark \ref{Rem:new}). Hence, $d$ has class number one or two. In the latter case, $LK=H(d)$ by class field theory. \qed

\begin{Remark}
For many K3 surfaces with Picard rank 20 over a quadratic extension, we know a model over $\Q$. Most of these models arise through the Shioda-Inose fibration (\cite{Inose}, \cite{S-fields}) or through extremal elliptic surfaces (\cite{BM}, \cite{S-Rocky}). It is an open question whether \emph{all} K3 surfaces with Picard rank 20 over a quadratic extension (or more generally with discriminant $d$ of class number two) might have a model over $\Q$.
\end{Remark}

\section{Singular K3 surfaces over number fields}
\label{s:gen}

We conclude the paper with an application of our techniques to general singular K3 surfaces. We will derive an explicit obstruction for the field of definition of the surface resp.~its N\'eron-Severi group. First we recall a possible field of definition:

\begin{Lemma}\label{Lem:def}
Let $X$ be a singular K3 surface of discriminant $d$. Then $X$ has a model over the ring class field $H(d)$.
\end{Lemma}

A model was given in \cite[proof of Proposition~10]{S-fields}, based on Inose's fibration in \cite{Inose} (cf.~\cite{Sandwich}). Elkies announced another model  in \cite{Elkies}.

In general, the field $H(d)$ need not be the optimal field of definition. In fact, there are examples of singular K3 surfaces over $\Q$ where $H(d)$ has degree 16 or 24 over $K$. The question arises how far one can possibly descend $X$, starting from $H(d)$. Shimada in \cite{Shimada} for fundamental $d$ and the author in \cite{S-fields} in full generality derived the following condition in terms of lattice: 
\begin{eqnarray}\label{eq:genus}
\{T(X^\sigma); \sigma\in\text{Aut}(\C/K)\} = \text{genus of } T(X).
\end{eqnarray}
In section~\ref{s:sing}, we used this to the following extent: If $X$ is defined over $\Q$, then the genus of $T(X)$ consists of a single class, i.e. $Cl(T(X))\cong (\Z/2)^g$. 

To rephrase (\ref{eq:genus}) in terms of class field theory, denote the degree of primitivity of $T(X)$ by $m$. Write $d=m^2 d'$, so that we can identify
\[
Cl(T(X))\cong Cl(d').
\]
Let $G=Cl(d')[2]$, the two-torsion subgroup of $Cl(d')$, and $M$ the fixed field of $G$ in the abelian Galois extension $H(d')/K$.

\begin{Lemma}\label{Lem:2}
Let $X$ be a singular K3 surface over some number field $L$. In the above notation, 
\[
M\subset KL.
\]
\end{Lemma}

So far, this was the only known obstruction to fields of definition of singular K3 surfaces. The only drawback of Lemma \ref{Lem:2} is that it fails to measure the degree of primitivity of $T(X)$. For this reason,  Theorem~\ref{thm:2} provides a major improvement: By providing bounds on the discriminant $d$, it also implies restrictions of the degree of primitivity. We shall now apply the techniques from the previous sections to prove  Theorem~\ref{thm:2}.

\emph{Proof of  Theorem~\ref{thm:2}:} Without loss of generality, we can assume that $L$ contains $K$. We consider all those good primes $p$ that split completely in $L$. Let $\p|p$ in $L$. Then $\F_\p=\F_p$ and $\rho(X/\F_p)=20$. Hence we can apply the Artin-Tate conjecture at $\p$. As in the previous section, Lemma \ref{Lem:shape} guarantees that $a_\p$ has the shape of Example~\ref{Ex:Hecke} and $p$ splits into prime ideals of order two in $Cl(K)$. By the argumentation of sections \ref{s:K}--\ref{s:d}, $p$ is not only represented by the principal class of $Cl(K)$, but also of $Cl(d)$. Hence, by class field theory, $L$ has to contain the ring class field $H(d)$. \qed

Since there are only limited possibilities for the Galois action on the N\'eron-Severi lattice of a singular K3 surface (or on any lattice of given rank),  Theorem~\ref{thm:2} provides us with a direct proof of the following finiteness result due to  \v Safarevi\v c. For best efficiency,  Theorem~\ref{thm:2} should be combined with Lemma~\ref{Lem:2}. 

\begin{Theorem}[\v Safarevi\v c {\cite{Shafa}}]\label{Thm:oo}
Let $n\in\N$. Then 
\[
\# \{\text{singular K3 surface $X$ over $L$}; \; [L:\Q]\leq n\}_{/\cong} < \infty.
\]
\end{Theorem}

\begin{Remark}
Similar results can be established for other modular surfaces, for instance for singular abelian surfaces (cf.~\cite{SM}). In that particular case, they would also follow from the cohomological structure (see Remark~\ref{Rem:ab}).
\end{Remark}


\vspace{0.5cm}

\textbf{Acknowledgement:} The author would like to thank N.~Elkies, I.~Shimada and T.~Shioda for many stimulating discussions. Funding from DFG under research grant Schu 2266/2-2 is gratefully acknowledged. 

\vspace{0.5cm}

\vspace{0.8cm}

Matthias Sch\"utt\\
Institute of Algebraic Geometry\\
Leibniz University Hannover\\
Welfengarten 1\\
30167 Hannover\\
Germany\\
schuett@math.uni-hannover.de
\end{document}